\documentclass[11pt]{amsart}
\usepackage{geometry}                
\usepackage{graphicx}
\usepackage{amssymb}
\usepackage{epstopdf}
\title{\textbf {Ideals on the Quantum Plane's jet space}}
\author{Andrey Glubokov}

\begin{document}

\maketitle

\textbf{The goal of this paper is to introduce some rings that play the role
of the jet spaces of the quantum plane and unlike the quantum plane itself
possess interesting nontrivial prime ideals. We will prove some results
(theorems 1-4) about the prime spectrum of these rings.}
\section{\textbf{Introduction}}

According to the classical perception of plane geometry the affine plane
corresponds to the algebra freely generated by two variables $x$ and $y$
subject to the trivial commutation relation $yx=xy.$ When the commutation
relation $yx=xy$ is replaced by $yx=qxy$ the resulting algebra is called the
quantum plane\cite{1},\cite{2}.

Objects like ''planes'' are expected to possess some analogue of ''curves''.
But the quantum plane possesses very few prime ideals. The idea of the paper
is to look at certain rings that play the role of jet spaces of the quantum
planes. This is done by introducing new ''jet'' variables in the style of
Kolchin's differential algebra\cite{3} and by considering commutation
relations among these variables which are compatible with the action of the
natural derivations on these rings. These are the multiplicative relations
unlike the ones of Weyl type considered in particular in \cite{4}. It turns
out these new rings possess plenty of prime ideals which are related to the
(commutative) geometry of $\Bbb{P}^{n}\times \Bbb{P}^{n},n\geq 1$
\section{\textbf{Background and Motivation.}}

\subsection{Quantum Symmetry (basic example)}

\smallskip It is well-known that symmetry plays an important role in modern
mathematics and theoretical physics. Gauge symmetry leads to the Standard
Model in high energy physics; crystallographic space symmetry is fundamental
to solid state physics, conformal symmetry is crucial to string theory and
critical phenomenon. In a sense, the progress of modern physics is
accompanied by study of symmetry.\bigskip

The mathematical counterpart of all the above mentioned symmetries and other
popular ones is group. Recently there has been great interest in the study
of quantum group and quantum symmetry. Quantum group in contrast to its
literal meaning is not a group, even not a semi-group. However a quantum
group is the deformation of the universal enveloping algebra of a
finite-dimensional semi-simple Lie algebra introduced by Drinfeld \cite{5}%
and Jimbo\cite{6} in their study of the Yang-Baxter equation.\bigskip 

The word ''quantum'' in quantum group is from the Yang-Baxter equation:
solutions of the classical Yang-Baxter equation are closely related with
classical or semi-simple groups, while solutions of the quantum Yang-Baxter
equation related with quantum groups. The quantum thus really differs from
the canonical quantization and possesses different meanings for different
systems. Conventionally, quantum group is a Hopf algebra which is neither
commutative nor cocommutative. A Hopf algebra is endowed with the algebra
homomorphisms: comultiplication $\Delta $ and counit $\varepsilon $, and the
algebra anti-homomorphism: the antipode $S$.

Quantum group theory has been developed in different directions. In the
quantum space approach \cite{1}, the initial object is a quadratic algebra
which is considered being as the polynomial algebra on a quantum linear
space. Quantum group appears like a group of automorphisms of the quantum
linear space.

The basic example is a Quantum Group $GL_{q}(2).$
Let $k$ be a ground field, $q\in k^{*}.$ By definition, the ring of
polynomial functions $F=F[GL_{q}(2)]$ is a Hopf algebra which can be
described in the following way. As a $k-$algebra, it is generated by $a,b,c,d
$ and a formal inverse of the central element
\[
D=DET_{q}\left( 
\begin{array}{ll}
a & b \\ 
c & d
\end{array}
\right) =ad-q^{-1}bc 
\]
where $a,b,c,d$ satisfy the following commutation relations:
\begin{eqnarray*}
ab &=&q^{-1}ba \\
ac &=&q^{-1}ca \\
cd &=&q^{-1}dc \\
bd &=&q^{-1}db \\
bc &=&cb \\
ad-da &=&(q^{-1}-q)bc
\end{eqnarray*}

The comultiplication $\Delta :F\rightarrow F\otimes F$ is defined by

\[
\Delta \left( 
\begin{array}{ll}
a & b \\ 
c & d
\end{array}
\right) =\left( 
\begin{array}{ll}
a & b \\ 
c & d
\end{array}
\right) \otimes \left( 
\begin{array}{ll}
a & b \\ 
c & d
\end{array}
\right) 
\]

where the tensor product in the r.h.s. denotes the usual product of matrices
in which products like $ab$ are replaced by $a\otimes b$ . The counit is
given by

\[
\varepsilon \left( 
\begin{array}{ll}
a & b \\ 
c & d
\end{array}
\right) =\left( 
\begin{array}{ll}
1 & 0 \\ 
0 & 1
\end{array}
\right) 
\]
The antipode map $S:F\rightarrow F$ is
\[
S\left( 
\begin{array}{ll}
a & b \\ 
c & d
\end{array}
\right) =D^{-1}\left( 
\begin{array}{ll}
d & -qb \\ 
-c/q & a
\end{array}
\right) 
\]
It can be checked directly that all these structures are well defined and
satisfy the Hopf algebra axioms.
\section{\textbf{Quantum plane: Gauss Polynomials and the q-Binomial Formula.}}

\textbf{Definition.} Let $k$ be a field . Let $q\neq 1$ be an invertible element of the ground
field $k$ and let $I_{q}$ be the two-sided ideal of the free algebra $k<x,y>$
of noncommutative polynomials in $x$ and $y$ generated by the element $%
f=yx-qxy.$ The quantum plane is defined as the quotient algebra $\frac{k<x,y>%
}{I_{q}}$ .
For future developments, we need to compute the powers of $x+y$ in the
quantum plane. To this end we have to consider Gauss polynomials..

Gauss polynomials are polynomials in one variable $q$ whose values at $q=1$
are equal to the classical binomial coefficients.
For any integer $n>0$ , set
\[
(n)_{q}=1+q+q^{2}+...+q^{n-1}=\frac{q^{n}-1}{q-1} 
\]
Define the $q-$factorial of $n$ by $(0)!_{q}=1$ and
\[
(n)!_{q}=(1)_{q}(2)_{q}...(n)_{q}=\frac{(q-1)(q^{2}-1)...(q^{n}-1)}{(q-1)^{n}
} 
\]

when $n>0.$ The $q-$factorial is a polynomial in $q$ with integral
coefficients and with value at $q=1$ equal to the usual factorial $n!.$ We
define the Gauss polynomials for $0\leq k\leq n$ by

\[
\binom{n}{k}_{q}=\frac{(n)!_{q}}{(k)!_{q}(n-k)!_{q}} 
\]

with following properties:

a). $\binom{n}{k}_{q}$is a polynomial in $q$ with integral coefficients and
with value at $q=1$ equal to the binomial coefficient $\binom{n}{k}.$

b).$q-$Pascal identity
\[
\binom{n}{k}_{q}=\binom{n-1}{k-1}_{q}+q^{k}\binom{n-1}{k}_{q}=\binom{n-1}{k}%
_{q}+q^{n-k}\binom{n-1}{k-1}_{q} 
\]
c). $q-$analog of the Chu-Vandermonde formula, for $m\geq p\geq n,$we
have 
\[
\binom{m+n}{p}_{q}=\sum\limits_{0\leq k\leq p}q^{(m-k)(p-k)}\binom{m}{k}_{q}%
\binom{n}{p-k}_{q} 
\]

Finally for all $n>0,$
\[
(x+y)^{n}=\sum\limits_{0\leq k\leq n}\binom{n}{k}_{q}x^{k}y^{n-k} 
\]

In particular, if $q$ is a root of unity of order $p>0,$
\[
(x+y)^{p}=x^{p}+y^{p} 
\]

Let $z$ be a variable commuting with $q.$The $q-$exponential can be defined
as following formal series: 
\[
e_{q}(z)=\sum\limits_{n\geq 0}\frac{z^{n}}{(n)!_{q}} 
\]
such that 
$e_{q}(x+y)=e_{q}(x)e_{q}(y)$

In this paper $q$ will eventually be assumed not a root of unity. However
some of the results can be extended to the case when $q$ is a root of unity
in which case the $q-$binomial formulae become relevant.

\section{\textbf{Quantum plane and Quantum Group.}}

A more conceptual approach to $GL_{q}(2)$ consists in introducing quantum
plane $\frac{k<x,y>}{I_{q}}$ and obtaining the commutation relations of $%
GL_{q}(2)$ from the following matrix relations:

\[
\binom{x^{\prime }}{y^{\prime }}=\left( 
\begin{array}{ll}
a & b \\ 
c & d
\end{array}
\right) \binom{x}{y} 
\]
\[
\binom{x^{\prime \prime }}{y^{^{\prime \prime }}}=\left( 
\begin{array}{ll}
a & c \\ 
b & d
\end{array}
\right) \binom{x}{y} 
\]

such that $x^{\prime },y^{\prime }$ and $x^{\prime \prime },y^{\prime \prime
}$ are on Quantum Plane and
\begin{eqnarray*}
y^{\prime }x^{\prime } &=&qx^{\prime }y^{\prime } \\
y^{\prime \prime }x^{\prime \prime } &=&qx^{\prime \prime }y^{\prime \prime }
\end{eqnarray*}

In this way, $GL_{q}(2)$ emerges merely as a quantum automorphism group of
noncommutative linear space.

\section{\textbf{The problem and the main results}}

The family of prime ideals of the quantum plane has a simple structure as we
shall presently review. Recall that an ideal $P$ is prime if $P\neq (1)$ and if from $ab\in P$ it follows that $a\in P$ or $b\in P.$ We denote by $Spec$ $B$ the set of prime ideals in any ring $B$.
$Spec\left( \frac{k<x,y>}{I_{q}}\right) $ consists of the following
prime ideals:
$\left\{ <0>,<x,y>,<x-\alpha ,y>,<x,y-\beta >\right\} $ where $\alpha ,\beta
\in k^{\times }$and $<S>$ denotes the two-sided ideal generated by set $S. $

Due to the commutation relation $yx=qxy$ the above set of ideals can be
rewritten as $\left\{ <0><x,y>,<x-\alpha >,<y-\beta >\right\} $ since, for
example, 
\[
y(x-\alpha )-q(x-\alpha )y=(q-1)\alpha y 
\]
so $y\in <x-\alpha >$ and $<x-\alpha ,y>=<x-\alpha >. $

The fact that the ring structure of the quantum plane is so trivial prevents
us from considering ''curves'' on it .
That is a motivation to attempt to introduce new rings that play the role of
the jet spaces of the quantum plane and possess interesting nontrivial prime
ideals.
Let us consider the noncommutative ring
$B^{(n)}=k<x,x^{\prime },x^{\prime \prime },...,x^{(n)},y,y^{\prime
},y^{\prime \prime },...y^{(n)}>$ where $x^{\prime },x^{\prime \prime
},...,x^{(n)},y^{\prime },y^{\prime \prime },...y^{(n)}$ are new
indeterminates.
Consider the unique $k-$derivation $\delta :$ $B^{(n-1)}\rightarrow B^{(n)}$
sending $\delta x=x^{\prime },\delta x^{\prime }=x^{\prime \prime },....$and 
$\delta y=y^{\prime },\delta y^{\prime }=y^{\prime \prime },...$
Recall that a $k-$derivation $\delta $ is a $k-$linear map satisfying the
usual Leibniz rule:
$\delta (FG)=\delta FG+F\delta G\
$

Define the following elements of $B^{(1)}:$

$g_{1}=y^{\prime }x-qxy^{\prime }$

$g_{2}=yx^{\prime }-qx^{\prime }y$

$g_{3}=y^{\prime }x^{\prime }-qx^{\prime }y^{\prime }$

Hence $\delta f=g_{1}+g_{2} $

Also define

$h=xx^{\prime }-x^{\prime }x, $
$\overline{h}=yy^{\prime }-y^{\prime }y$
and, more generally, define the following elements of $B^{(n)}:$

$h_{ij}=x^{(i)}x^{(j)}-x^{(j)}x^{(i)}$

$\overline{h_{ij}}=y^{(i)}y^{(j)}-y^{(j)}y^{(i)}$

$g_{ij}=y^{(i)}x^{(j)}-qx^{(j)}y^{(i)}$for integers $i,j\leq n.\medskip $
 Finally, define $A^{(n)}=\frac{B^{(n)}}{<f,h_{ij},\overline{h_{ij}}%
,g_{ij}>}$

For further purposes we define the commutative ring
$A_{c}^{(n)}=k[x,x^{\prime },x^{\prime \prime },...,x^{(n)},y,y^{\prime
},y^{\prime \prime },...y^{(n)}]$
which is the ring of the usual(commutative) polynomials.
Any monomial $x^{i_{o}}$ $(x^{\prime })^{i_{1}}$ $(x^{\prime \prime
})^{i_{2}}$ $...$ $(x^{(n)})^{i_{n}}$ $y^{j_{0}}(y^{\prime })^{j_{1}}$ $%
(y^{\prime \prime })^{j_{2}}$ $...(y^{(n)})^{j_{n}}$has the \textbf{%
bi-degree }$(i,j)$ where the total degree in $x,x^{\prime },...,x^{(n)}$ is $%
i=i_{0}+i_{1}+...+i_{n}$ and the total degree in $y,y^{\prime },y^{\prime
\prime },...y^{(n)}$is $j=j_{0}+j_{1}+...+j_{n}$

\textbf{Lemma.} There is a unique $k-$linear bijective map $A^{(n)}\cong
A_{c}^{(n)}$sending the class of any monomial into the same monomial viewed
as an element of $A_{c}^{(n)}.$

Via the above bijection we have a multiplication law $\bullet _{c}$ on $%
A_{c}^{(n)}$such that for any two bi-homogeneous polynomials of bi-degrees $
(i,j)$ and $(k,l)$ respectively,
\[
a\bullet _{c}b=ab\cdot q^{-jk}. 
\]

The bijection in the Lemma is not an isomorphism of rings.
From now on we shall identify $A_{c}^{(n)}$and $A^{(n)}$as sets via above
bijection. Note that $A_{c}^{(n)}$is bi-graded in the usual way. In the following let $q$ be not a root of unity. Our main results about $Spec\left( A^{(n)}\right) $ can be presented as the
following theorems 1-4.\medskip 

\textbf{Theorem 1. }

If $0\neq P\subset A^{(n)}$\ is a prime ideal then $P$ contains a non-zero
bi-homogeneous polynomial which as an element of $A_{c}^{(n)}$\ is
irreducible.

\smallskip

\textbf{Theorem 2.}

If $f\in A^{(n)}$ bi-homogeneous such that its image $f\in A_{c}^{(n)}$ is
irreducible

then $<f>\subset A^{(n)}$ is prime.

\smallskip

\textbf{Theorem 3. }Any prime ideal $P\subset A^{(n)}$ not containing any of the ideals $%
<x,x^{\prime },...,x^{(n)}>$ or $<y,y^{\prime },....y^{(n)}>$ is of the form 
$P=<T>,$ where $T$ is the family of all bi-homogeneous polynomials in $P.$

\textbf{Theorem 4.a }Any prime ideal $P\subset A^{(n)}$  such that $<x,x^{\prime },...,x^{(n)}>$ $
\subset P$
is of the form $<x,x^{\prime },...,x^{(n)},\varphi _{1}(y,y^{\prime
},....y^{(n)}),...,\varphi _{k}(y,y^{\prime },....y^{(n)})>$
where $\varphi _{i}(y,y^{\prime },....y^{(n)})\subset k[y,y^{\prime
},....y^{(n)}]$ for $i=1,...,k$ generate a prime ideal of $k[y,y^{\prime
},....y^{(n)}].$

\textbf{Theorem 4.b}Any prime ideal $P\subset A^{(n)}$  such that $<y,y^{\prime },....y^{(n)}>$ $%
\subset P$
is of the form $<\psi _{1}(x,x^{\prime },...,x^{(n)}),...\psi
_{k}(x,x^{\prime },...,x^{(n)}),y,y^{\prime },....y^{(n)})>$
where $\psi _{i}(x,x^{\prime },...,x^{(n)})\subset k[x,x^{\prime
},...,x^{(n)}]$ for $i=1,...,k$ generate a prime ideal of $k[x,x^{\prime
},...,x^{(n)}].$

\section{$\delta -$\textbf{prime ideals}}

Let us recall the previously defined derivation $\delta
:B^{(n-1)}\rightarrow B^{(n)}.$ Let $A=\underrightarrow{\lim }$ $A^{(n)}$.
Then $\delta $ induces a $k-$derivation $\delta :A\rightarrow A.$
For each $n$ we have $x^{(n)}=\delta x^{(n-1)}$ and $y^{(n)}=\delta
y^{(n-1)}$

Define a $\delta -$prime ideal to be a prime ideal $P$ such that $\delta
P\subset P.$ As in Theorem 3 let $T=\left\{ f\in P\left| f\text{ is bi-homogeneous}%
\right. \right\} $ so $P=<T>$

We can prove the following \textbf{Proposition}:

$$\delta P\subset P\Longleftrightarrow \delta T\subset T$$

Proof.

Implication $\Longrightarrow :$ $\delta T\subset \delta P\cap \{$bi-homogenious
elements of $A\}\subset P\cap \{$bi-homogenious elements of $A\}=T$

The $\Longleftarrow $ part follows because if $f\in P$ then

$$f=\sum_{i}\alpha _{i}f_{i}\beta _{i},f_{i}\in T$$ 
so $$\delta f=\sum _{i}\delta \alpha _{i}f_{i}\beta _{i}+\sum _{i} \alpha _{i}\delta f_{i}\beta _{i}+\sum _{i} \alpha _{i}f_{i}\delta \beta _{i}$$

 Since $\delta T\subset T$ , then $\delta f\in <T>=P$

\section{Appendix}
For the proofs of theorems 1-4 we need the following definition of the lexicographical
ordering\textbf{\ }in $\Bbb{N}^{2}$

\[
(a,b)\leq (c,d)\Longleftrightarrow \left\{ 
\begin{array}{c}
\text{either }a<c \\ 
\text{or }a=c\text{ },b\leq d
\end{array}
\right\} 
\]
\medskip

Let's consider a polynomial $g\in A^{(n)}$
Write $g=\sum_{ij}g_{ij}$, such that $g_{ij}$ is bi-homogeneous of bi-degree $%
(i,j).$ Let's consider the set $\Gamma _{g}=\left\{ (i,j),g_{ij}\neq 0\right\}
\subset \Bbb{N}^{2}$.
The \textbf{size }of a polynomial $g$
in  $A^{(n)}$ is size(g)= $\#\Gamma _{g}$ , number of
points in $\Gamma _{g}.$ If $g_{ij}$ has the bi-degree $(i,j)$ then $g_{ij}x$
has the bi-degree $(i+1,j)$ and $yg_{ij}$ has bi-degree $(i,j+1).$ The size
of $gx$ and $yg$ will stay the same as the size of $g_{{}}.$

\textbf{Lemma 1. }If $h=yg-q^{\nu }gy$ and $(\nu ,\mu )\in $ $\Gamma _{g}$ ,
then the size of $\Gamma _{h}$ will be strictly less than the size of $%
\Gamma _{g}.$

Indeed, $$h=yg-q^{\nu }gy=y \sum _{ij}g_{ij}-q^{\nu } \sum _{ij}g_{ij}y=\sum _{ij}g_{ij}y(q^{i}-q^{\nu })$$ . 

It follows that all points of $\Gamma _{g}$with the first coordinate equal to $\nu $
will disappear in $\Gamma _{h}$ and the size of $\Gamma _{h}$ will be
strictly less than the size of $\Gamma _{g}.\medskip $

Similarly, if $\overline{h}=gx-q^{\mu }xg$ and $(\nu ,\mu )\in $ $\Gamma
_{g}, $ then the size of $\Gamma _{\overline{h}}$ will be strictly less than
the size of $\Gamma _{g}.$
\smallskip
\subsection{\textit{Proof of Theorem 1. }}

We start by showing the following claim: there exists a nonzero
bi-homogeneous polynomial in\textbf{\ }$P.$
Indeed take $0\neq g\in P$ of smallest possible size. We claim that $%
size(g)=1$ which means $g$ is bi-homogeneous.
Assume that $size(g)\geq 2$

Case 1. $g$ is not homogeneous in $x,x^{\prime },...,x^{(n)}.$

Let's consider $\overline{g}=yg-q^{\nu }gy\in P$ such that there is at least
one term with total degree in $x,x^{\prime },...,x^{(n)}$ equal to $\nu .$
Since $g$ is not homogeneous in $x,x^{\prime },...,x^{(n)},
\overline{g}\neq 0.$ On the other hand by the Lemma we have $\#\Gamma _{%
\overline{g}}<\#\Gamma _{g}$ which contradicts the minimality of $size(g).$

Case 2. $g$ is homogeneous in $x,x^{\prime },...,x^{(n)}$ but not in $%
y,y^{\prime },....y^{(n)}.$

Let's consider $\hat{g}=gx-q^{\mu }xg\in P$ such that there is at least one
term in $g$ with the total degree in $y,y^{\prime },....y^{(m)}$ equal to $%
\mu .$
Since $g$ is not homogeneous in $y,y^{\prime },....y^{(n)}$, $\hat{g}\neq 0. 
$ On the other hand by the Lemma we have $\#\Gamma _{\hat{g}}<\#\Gamma _{g}$
which contradicts with minimality of $size(g)$

This proves our claim. To conclude the proof of the Theorem 1, using our
claim one can pick a nonzero bi-homogeneous polynomial $f\in P$ of smallest
bi-degree $(i^{*},j^{*})$ with respect to lexicographical order among the
nonzero bi-homogeneous polynomials in $P$.

We claim that $f$ is irreducible in $A_{c}^{(n)}.$

If we assume it is not,
then $f=g\cdot _{c}h,$ $g,h\in A_{c}^{(n)},($ $g,h\notin k)$

Write:
\begin{eqnarray*}
g &=&\sum_{ij}g_{ij}\text{ ,}g_{ij}\text{ bi-homogeneous of bidegree }(i,j)
\\
h &=&\sum_{ij}h_{ij}\text{, }h_{ij}\text{ bi-homogeneous of bidegree }(i,j)
\end{eqnarray*}

Note the following properties of bi-degrees:\medskip

1. bideg$(g_{i_{1}j_{1}}\cdot h_{i_{2}j_{2}})=(i_{1}+i_{2},j_{1}+j_{2})$

\smallskip
2. 
\begin{eqnarray*}
\text{If }(i_{1},j_{1}) &<&(i_{0},j_{0}) \\
\text{and }(k_{1},l_{1}) &\leq &(k_{0},l_{0}) \\
\text{then }(i_{1}+k_{1},j_{1}+l_{1}) &<&(i_{0}+k_{0},j_{0}+l_{0})
\end{eqnarray*}

Let $(i_{0},j_{0})$ be the highest element of $\Gamma _{g}$ with respect to
lexicographical order$,(k_{0},l_{0})$ be the highest element of $\Gamma _{h}$
with respect to lexicographical order$.$
and let $(i_{1},j_{1})$ be the lowest element of $\Gamma _{g},(k_{1},l_{1})$ be the lowest element of $\Gamma _{h}\ $

Then the highest element of $\Gamma _{g\text{ }\cdot _{c}h}$ will be $%
(i_{0}+k_{0},j_{0}+l_{0})$ and the lowest element of $\Gamma _{g\text{ }%
\cdot _{c}h}$ will be $(i_{1}+k_{1},j_{1}+l_{1}).$ Since $f=g\cdot _{c}h$ we
have $(i_{0}+k_{0},j_{0}+l_{0})=(i_{1}+k_{1},j_{1}+l_{1})=(i^{*},j^{*})%
 $

Since $i^{*}=i_{0}+k_{0}=i_{1}+k_{1}$and $i_{0}\geq i_{1}$ it follows that $%
i_{0}=i_{1}$ because if $i_{0}>i_{1}$ then $k_{0}$ has to be \textbf{less }%
then $k_{1}$ which contradicts with the choice of $k_{0}.$ It immediately
follows that $k_{0}=k_{1}. $
Similarly, $j_{0}=j_{1}$ and $l_{0}=l_{1}$, so $g$ and $h$ are both bi-homogeneous of degrees less than $%
(i^{*},j^{*}).\medskip $

Since $P$ is a prime ideal, at least one of them belongs to $P.$ This
contradicts the choice of $f.$

\[
\]

\subsection{\textit{Proof of the Theorem 2.\protect\medskip }}

Assume $f$ is irreducible in $A_{c}^{(n)}$ and bi-homogeneous of bi-degree $
(i,j). $

We prove by induction on the total degree $N$ in $x,x^{\prime },...,x^{(n)},$
$y,y^{\prime },....y^{(n)}$ that if $f$ has a total degree $N$ then from $%
g\cdot h\in <f>$it follows that $g$ or $h\in <f> $

If $N=0$ the theorem is clear.
Assume the theorem is true for total degree less or equal to $N-1.$

Let $N$ be the total degree of $f.$
We have that from $g\cdot h\in <f>$it follows that $g$ $\cdot $ $%
h=\sum_{i}\alpha _{i}f\beta _{i}$ where $\alpha _{i\text{ }}$ and $\beta
_{i} $ belong to $A^{(n)}.$We may assume that $\alpha _{i\text{ }}$and $%
\beta _{i} $ are bi-homogeneous.\medskip

Since $f$ is bi-homogeneous, $\sum_{i}\alpha _{i}f\beta _{i}=$ $%
\sum_{i}\alpha _{i}\beta _{i}q^{n_{i}}f=\sum_{i}\gamma _{i}f=\gamma f,$ for
some $n_{i},\gamma _{i},\gamma .$

\begin{eqnarray*}
\gamma &=&\sum_{ij}\gamma _{ij}\text{, }\gamma _{ij}\text{ bi-homogeneous of
bidegree }(i,j)
\end{eqnarray*}

Let $(i_{0},j_{0})$ be the highest element of $\Gamma _{g}$ with respect to
lexicographical order$,(k_{0},l_{0})$ be the highest element of $\Gamma _{h}$
and $(m_{0},n_{0})$ be the highest element of $\Gamma _{\gamma }.$ Then
\begin{eqnarray*}
g_{i_{0}j_{0}}\cdot h_{k_{0}l_{0}} &=&\gamma _{m_{0}n_{0}}\cdot f \\
q^{t}g_{i_{0}j_{0}}\cdot _{c}h_{k_{0}l_{0}} &=&q^{s}\gamma
_{m_{0}n_{0}}\cdot _{c}f
\end{eqnarray*}

for some $t$ and $s.$

Since $f$ is irreducible in the commutative ring $A_{c}^{(n)},$ it follows
that $g_{i_{0}j_{0}}=\eta \cdot _{c}f=q^{-li}\eta \cdot f$ ($f$ is
bi-homogeneous and the bi-degree of $\eta $ is $\left( k,l\right) $ ) or $%
h_{k_{0}l_{0}}=\overline{\eta }\cdot _{c}f=q^{-\overline{l}i}\overline{\eta }%
\cdot f$ $($bi-degree of $\overline{\eta \text{ }}$is $\left( \overline{k},%
\overline{l}\right) \medskip )$

Assume, for example, the former is the case.
From $gh=\gamma f$ we get $(g-g_{i_{0}j_{0}}+g_{i_{0}j_{0}})\cdot
h=(g-g_{i_{0}j_{0}})\cdot h+g_{i_{0}j_{0}}\cdot h=g^{\prime }\cdot
h+q^{w}\eta \cdot f\cdot h$ where $g^{\prime }=g-g_{i_{0}j_{0}}.$Obviously, $%
g^{\prime }\cdot h\in <f>.\medskip $

Since the total degree in $x,x^{\prime },...,x^{(n)},$ $y,y^{\prime
},....y^{(m)}$of $g^{\prime }\cdot h$ is less or equal to $N-1$, by the
induction hypothesis either $g^{\prime }$ $\in <f>$and $g=g^{\prime
}+q^{w}\eta \cdot f$ $\in <f>$or $h\in <f>$ and we are done.

\subsection{\textit{Proof of the theorem 3. }}

It is obvious that $<T>\subset P.$

To prove $<T>\supset P$ assume on the contrary that $P$ does not belong to $
<T>.$
Let $f\in P\setminus <T>$ be of minimal size. Since by this assumption $f$
cannot be bi-homogeneous the $size(f)>1.$ There are two cases.\medskip

\textbf{Case1.}
$f$ is not homogeneous in $x,x^{\prime },...,x^{(n)}$.
Write $f=\sum_{st}f_{st}.$
For an arbitrary $(k,l)\in \Gamma _{f}.$

There exists a pair $(i,j)\in \Gamma _{f}$ such that $i\neq k$ otherwise $f$
ought to be homogeneous in $x,x^{\prime },...,x^{(n)}.\medskip $

Let $h=yf-q^{i}fy\in P.$

Then by the Lemma 1 $size(h)<size(f).$ 
It follows that $h\in <T>$ so $h$ can be written as
\begin{eqnarray*}
h &=&\gamma _{1}B_{1}+\gamma _{2}B_{2}+...+\gamma _{m}B_{m} \\
\text{where }B_{l}\text{ } &\in &P\text{ and bi-homogeneous }
\end{eqnarray*}

Also $h=\sum_{st}(yf_{st}-q^{i}f_{st}y)=\sum_{st}\lambda
_{si}f_{st}y$, where $\lambda _{si}=q^{s}-q^{i}.$

\medskip

Let us pick out the bi-homogeneous components of bi-degree $(k,l+1).$
Then $\lambda _{ki}$ $\cdot f_{kl}\cdot y=\widetilde{\gamma _{1}}B_{1}+%
\widetilde{\gamma _{2}}B_{2}+...+\widetilde{\gamma _{m}}B_{m}$ $\in <T>$
where $\widetilde{\gamma _{1}},\widetilde{\gamma _{2,}}....,\widetilde{
\gamma _{m}}$ are bi-homogeneous. 
Since $\lambda _{ki}\neq 0$ because $i\neq k,$ we have
$
f_{kl}\cdot y\in <T>. 
$
So $f_{kl}\cdot y\in P\medskip $

Similarly let $h^{(s)}=y^{(s)}f-q^{i}fy^{(s)}$ .
As above we get $f_{kl}\cdot y^{(s)}\in P$ for all $s$.
Since $<y,y^{\prime },....y^{(n)}>$ is not contained in $P$ it follows that
at least one of
$y^{(s)}$ $\notin P.$
Because $P$ is prime, $f_{kl}\in P.$ But $f_{kl}$ is obviously bi-homogeneous so $f_{kl}\in <T>$

Since the pair $(k,l)$ is arbitrary it follows that $f=\sum
_{st}f_{st}\in <T>$ - a contradiction.\medskip

\textbf{Case 2. }
$f$ is homogeneous in $x,x^{\prime },...,x^{(n)}$ but not in $y,y^{\prime
},....y^{(n)}.\medskip $
Write $$f=\sum_{t}f_{kt}$$
For an arbitrary $(k,l)\in \Gamma _{f}.\medskip $
there exists a pair $(k,j)\in \Gamma _{f}$ such that $l\neq j$ otherwise $f$
ought to be homogeneous in $y,y^{\prime },....y^{(n)}.\medskip $

Let $\overline{h}=fx-q^{j}xf\in P.$
Then by the Lemma 1 $size(\overline{h})<size(f).$ 
It follows that $\overline{h}\in <T>$ so $\overline{h}$ can be written as
\begin{eqnarray*}
\overline{h} &=&\overline{\gamma }_{1}\overline{B}_{1}+\overline{\gamma }_{2}%
\overline{B}_{2}+...+\overline{\gamma }_{m}\overline{B}_{m} \\
\text{where }\overline{B_{l}}\text{ } &\in &P\text{ and bi-homogeneous .}
\end{eqnarray*}
Also $\overline{h}=\sum_{t}(f_{kt}x-q^{j}xf_{kt})=\sum_{t}
\lambda _{tj}f_{kt}y$, where $\lambda _{si}=q^{t}-q^{j}.$
Let us pick out the bi-homogeneous components of bi-degree $(k+1,l). 
$
Then $\lambda _{lj}$ $\cdot x\cdot f_{kl}=\overline{\widetilde{\gamma _{1}}}%
\overline{B_{1}}+\overline{\widetilde{\gamma _{2}}}\overline{B_{2}}+...+%
\overline{\widetilde{\gamma _{m}}}\overline{B_{m}}$ $\in <T>$ where $%
\overline{\widetilde{\gamma _{1}}},\overline{\widetilde{\gamma _{2,}}}....,%
\overline{\widetilde{\gamma _{m}}}$ are bi-homogeneous.

Since $\lambda _{lj}\neq 0$ because of$l\neq j,$ we have
$x\cdot f_{kl}\in <T> $,
so $x\cdot f_{kl}\in P$

Similarly let $\overline{h^{(s)}}=fx^{(s)}-q^{j}x^{(s)}f$ .
As above we get $x^{(s)}\cdot f_{kl}\in P$ for all $s$.
Since $<x,x^{\prime },...,x^{(n)}>$ is not contained in $P$ it follows that
at least one of
$x^{(s)}$ $\notin P. $
Since $P$ is prime, $f_{kl}\in P, $
but $f_{kl}$ is obviously bi-homogeneous so $f_{kl}\in <T>. $
The pair $(k,l)$ is arbitrary so it follows that $f=\sum_{st}f_{st}\in <T>$ - a contradiction.

\subsection{\textit{Proof. of Theorem 4.a}}

Let us consider the factor ideal $\frac{P}{<x,x^{\prime },...,x^{(n)}>}.$
Then
\[
\frac{P}{<x,x^{\prime },...,x^{(n)}>}\subset \frac{k<x,x^{\prime
},...,x^{(n)},y,y^{\prime },....y^{(n)}>}{<x,x^{\prime },...,x^{(n)}>}%
=k[y,y^{\prime },....y^{(n)}] 
\]

Due to the structure of prime ideals of $k[y,y^{\prime },....y^{(n)}]$ we
have

$\frac{P}{<x,x^{\prime },...,x^{(n)}>}=\left( \varphi _{1}(y,y^{\prime
},....y^{(n)}),...,\varphi _{k}(y,y^{\prime },....y^{(n)})\right) \cdot 
\frac{k<x,x^{\prime },y,y^{\prime },....y^{(m)}>}{<x,x^{\prime },...,x^{(n)}>%
}=\medskip $

$\frac{<x,x^{\prime },...,x^{(n)},\varphi _{1}(y,y^{\prime
},....y^{(n)}),...,\varphi _{k}(y,y^{\prime },....y^{(n)}>}{<x,x^{\prime
},...,x^{(n)}>}\medskip $

It follows that 
\[
P=<x,x^{\prime },...,x^{(n)},\varphi _{1}(y,y^{\prime
},....y^{(n)}),...,\varphi _{k}(y,y^{\prime },....y^{(n)}> 
\]

\text{Theorem 4.b can be proved similarly.}
\[
\]

\[
\]
\maketitle{   Ave Maria University, FL, USA, e-mail: andrey.glubokov@avemaria.edu}


\begin{thebibliography}{9}
\bibitem{1}  Yu. I. Manin.Quantum groups and noncommutative geometry.
Montreal University, 1988

\bibitem{2} Yu. I. Manin, Notes on quantum groups and quantum de Rham
complexes, Teoret. Mat. Fiz. 92 (1992), no. 3, 425--450; translation in
Theoret. and Math. Phys. 92 (1992), no. 3, 997--1023 (1993)

\bibitem{3} E. R. Kolchin, Differential Algebra and Algebraic Groups,
Academic Press, New York, 1973

\bibitem{4} Artin, M.; Tate, J.; Van der Bergh, M., Modules over
regular algebras of dimention 3, Inventiones mathematicae, volume 106;
pp.335-338; 1991

\bibitem{5}  V. G. Drinfeld, in Proceedings of the International
Congress of Mathematicians,Berkeley 1986, ed. A. M. Gleason (American
Mathematical Society, Providence,RI, 1987) p798.

\bibitem{6}  M. Jimbo, Lett. Math. Phys. 10 (1985) 63.
\end{thebibliography}
\end{document}